\documentclass[12pt]{article}
\usepackage{amssymb, amsmath}

\newtheorem{thm}     {Theorem}[section]
\newtheorem{prop}    [thm]{Proposition}

\newtheorem{cor}     [thm]{Corollary}
\newtheorem{lemma}   [thm]{Lemma}
\newtheorem{remark}   [thm]{Remark}
\newtheorem{example}   [thm]{Example}

\newcommand{\proof} {\noindent{\bf Proof. }}
\newcommand{\C}{\mathbb C}
\newcommand{\D}{\mathbb D}
\newcommand{\R}{\mathbb R}
\def\Re{{\rm Re\,}}
\def\bar{\overline}
\def\dist{{\rm dist}}

\def\const{{\rm const}}
\def\({\left(}
\def\){\right )}
\def\pv{{\rm p.v.}}
\def\id{{\rm id}}

\begin{document}

\title{Commutators of singular integrals, \\
the Bergman projection,
and boundary\\ regularity of elliptic equations in the plane}
\author{Alexander Tumanov
\footnote{The author is partially supported by a grant from the Simons Foundation.}}
\date{}
\maketitle

{\small
University of Illinois, Department of Mathematics,
1409 West Green Street, Urbana, IL 61801, USA,
tumanov@illinois.edu
}
\medskip

Abstract.
We obtain estimates of commutators of singular
integral operators in Lipschitz spaces and apply
the results to boundary regularity of elliptic
equations in the plane. We obtain an explicit asymptotic
formula for the Bergman projection.
\medskip


\medskip

\section{Introduction}

In this paper we are concerned with sharp boundary
regularity in Lipschitz spaces $C^{k,\alpha}$ of
first order elliptic equations of the form
\begin{equation}
\label{40}
f_{\bar z}=a(z)f_z+b(z)\bar f_{\bar z}+c(z)
\end{equation}
in a smooth bounded domain $\Omega$ in complex plane $\C$.
Here $f_z=\partial f/\partial z$ and
$f_{\bar z}=\partial f/\partial\bar z$.
We impose the ellipticity condition
\begin{equation}
\label{42}
|a(z)|+|b(z)|\le a_0<1
\end{equation}
for some constant $a_0$.
We first consider the scalar equation \eqref{40}
with the Dirichlet type boundary condition
$\Re f|_{b\Omega}=f_0$ for given function $f_0$
on the boundary. Suppose $a, b, c$ are in
$C^{k,\alpha}(\Omega)$, $k\ge 0$, $0<\alpha<1$.
Suppose $f_0$ is in $C^{k+1,\alpha}(b\Omega)$.
We would like to conclude that every generalized
solution of \eqref{40} with $\Re f|_{b\Omega}=f_0$
is automatically in $C^{k+1,\alpha}(\Omega)$.
Apparently, this classical question is not covered
in the extensive literature on the subject.
In particular, the case of first order equations
does not follow from the classical results on boundary
regularity of elliptic equations \cite{ADN, S}.
If in the scalar equation \eqref{40}, the coefficient
$b=0$, then the conclusion is rather simple
(see \cite{ST}, Proposition 2.1).
Tadeusz Iwaniec explained to the author that
for $k\ge 1$ the equation \eqref{40} can be reduced
to a second order equation, and then the conclusion follows
from Schauder's theory \cite{S}. In this paper
we give a proof for all $k\ge 0$.
We also consider the vector version of \eqref{40}
and give a proof of the regularity of the Dirichlet
problem for $a=0$.

Another common boundary condition for equation \eqref{40}
is $K_\Omega f=f_0$. Here $K_\Omega$ is the Cauchy
type integral \eqref{41} and $f_0$ is a given holomorphic
function in $\Omega$. In particular, the homogeneous
condition $K_\Omega f=0$ means that $f$ holomorphically
extends to $\C\setminus\Omega$ and vanishes at infinity.
Solving \eqref{40} with this boundary condition is
equivalent to the problem of inverting the operator
$f\mapsto f-T_\Omega(af_z+b\bar f_{\bar z})$, here
$T_\Omega$ is the Cauchy--Green operator \eqref{17}.
The vector version of this problem
with $a=0$ and small $b$ arises in
constructing small pseudoholomorphic curves
(see \cite{BGR, IvRo, NiWo}).
We prove the boundary regularity of this
problem in the scalar case for general $a$ and $b$
satisfying \eqref{42} and in the vector case for
$a=0$ and $\|b\|_\infty<1$, answering a question raised
in \cite{BGR}.

A classical approach \cite{AIM, V} to equation \eqref{40},
in particular, the Beltrami equation consists of reducing
\eqref{40} to an integral equation with the operator
$S_\Omega$ given by \eqref{43} or its modifications.
The solution operator of the integral equation is bounded
in $L^p(\Omega)$ for $p$ close to 2. In this approach,
it is essential that $\|S_\Omega\|_2\le 1$.
However, there is more precise information about the operator
$S_\Omega$, in particular, $S_\C$ is an isometry
of $L^2(\C)$, that is, $S_\C\bar S_\C=I$.
There is a related property of $S_\Omega\bar S_\Omega$
that we derive in Section 4.
An iteration of the integral equation corresponding to
\eqref{40} with $b\ne0$ involves the term
$S_\Omega b\bar S_\Omega\bar b$, here $b$ denotes the
operator of multiplication by $b$.
Since $S_\Omega$ and $\bar S_\Omega$ do not stand next
to each other, in order to make use of
$S_\Omega\bar S_\Omega$ we need information
about the commutators of $S_\Omega$ with multiplication
operators, namely, their smoothing properties.

There are well known $L^p$ estimates of commutators
of singular integral operators with multiplication
operators (see, e. g., \cite{C, CRW}).
However, apparently, $C^{k,\alpha}$ estimates of the
commutators are covered in the literature only for
the case of Cauchy type integrals and similar operators
(see \cite{P}, Section 3.4.1).
We present results on the matter for the operator $S_\Omega$.
Although we use complex variable notations,
the results are real in nature and could be established for
more general Calder\'on--Zigmund operators.

As we mentioned above, the scalar equation \eqref{40} with
$b=0$ is rather simple. If $b=0$, then the equation \eqref{40}
can be reduced to the case $a=b=0$ by changing the independent
variable. Our method involving $S_\Omega\bar S_\Omega$ and
the commutators now lets us deal with the case $a=0$, $b\ne0$.
In the scalar case, it suffices for treating the general
equation \eqref{40} because we can reduce it to $a=0$.
However, in the vector case, obviously, the reduction
to $a=0$ by changing the independent variable is not possible
in general, thus we only handle
the vector equation \eqref{40} for $a=0$.

In Sections 2 and 3 we include results on $C^{k,\alpha}$
regularity of commutators of $S_\Omega$ with multiplication
operators.
In Section 4 we study properties of $S_\Omega\bar S_\Omega$.
In Section 5 we give an asymptotic formula of the Bergman
projection for $\Omega$ in terms of $S_\Omega\bar S_\Omega$.
In Section 6
we treat integral equations corresponding to
the vector version of the equation \eqref{40} with $a=0$.
Finally, in Sections 7 and 8 we study the boundary
regularity of the equation \eqref{40}.

I wish to thank Tadeusz Iwaniec for his letter with
a sketch of the proof of Theorem \ref{Dirichlet-scalar}
for the case $k=0$, $b=0$ and the case $k\ge1$.
I am also grateful to Steve Bell for discussions
on the Bergman projection.
Finally, I thank Elias Stein for answering my
inquiry regarding singular integrals.

\section{Commutators of singular integrals}

For a domain $\Omega\subset\C$, we consider
the Calder\'on--Zigmund operator (see \cite{AIM, V})
\begin{equation}
\label{43}
S_\Omega u(z)=\pv\int_\Omega \frac{u(t)\,d^2t}{(t-z)^2}.
\end{equation}
Here for brevity $d^2t=(2\pi i)^{-1}dt\wedge d\bar t$,
and the integral is understood as Cauchy principal value.
Let $a(z)$ be a function in $\Omega$. We use the same notation
$a$ for the operator of multiplication by $a$.
We are concerned with smoothing properties of the commutator
\begin{equation}
\label{16}
[S_\Omega,a]u(z)=\int_\Omega \frac{a(t)-a(z)}{(t-z)^2}\,
u(t)\,d^2t
\end{equation}
in Lipschitz spaces. As usual $C^{k,\alpha}(\Omega)$ denotes
the space of functions whose derivatives to order $k\ge0$
satisfy a Lipschitz condition with exponent $0<\alpha<1$.
We also sometimes write $C^\alpha(\Omega)=C^{0,\alpha}(\Omega)$
and $C^{k+\alpha}(\Omega)=C^{k,\alpha}(\Omega)$.
We do not make a difference between
$C^{k,\alpha}(\Omega)$ and $C^{k,\alpha}(\bar\Omega)$.
If $k$ is integer, then we use $C^{k}(\Omega)$ for
the usual $C^k$-smooth functions.
We use $C^{k,\alpha}(\Omega,\R)$ for the set of real valued
functions in $C^{k,\alpha}(\Omega)$.
We denote by $L^p_R$ and $C^{k,\alpha}_R$ the spaces
of functions respectively in $L^p(\C)$ and $C^{k,\alpha}(\C)$
with support in the disc $|z|\le R$.
\begin{thm}
\label{RegSC}
Let $0<\alpha<1$, $0<\beta<\alpha$, $R>0$.
Let $S=S_\C$.
\begin{itemize}
\item[(i)]
If $a\in C^\alpha(\C)$,
then the commutator $[S,a]$ is a bounded operator
$L^\infty_R\to C^\beta(\C)$ and
$C^\beta_R\to C^\alpha(\C)$.
\item[(ii)]
If $a\in C^{k+1,\alpha}(\C)$, $k\ge 0$,
then $[S,a]$ is a bounded operator
$C^{k,\alpha}_R \to C^{k+1,\alpha}(\C)$.
\end{itemize}
\end{thm}

We begin the proof with a simple formula. We introduce
the difference and shift operators
\[
\Delta_h f(z)=f(z+h)-f(z),
\qquad
\delta_h f(z)=f(z+h).
\]
\begin{lemma}
\label{Delta-Commutator}
For the operator $S=S_\C$, the following formula holds.
\begin{equation}
\label{3}
\Delta_h [S,a]u=[S,\Delta_h a]\,\delta_h u
+[S,a]\,\Delta_h u.
\end{equation}
\end{lemma}
\proof
Clearly
$[\Delta_h,S]=0$, $[\delta_h,S]=0$, and
$[\Delta_h,a]=(\Delta_h a)\delta_h$.
Then
\begin{align*}
\Delta_h [S,a]
&=\Delta_h Sa-\Delta_h aS=S\Delta_h a-\Delta_h aS
\\
&=S([\Delta_h,a]+a\Delta_h)-([\Delta_h,a]+a\Delta_h)S
\\
&=S((\Delta_h a)\delta_h+a\Delta_h)-((\Delta_h a)\delta_h+a\Delta_h)S
\\
&=[S(\Delta_ha)\delta_h-(\Delta_ha)S\delta_h]
+[Sa\Delta_h-aS\Delta_h]
=[S,\Delta_h a]\,\delta_h +[S,a]\,\Delta_h.
\quad
\blacksquare
\end{align*}

\noindent
{\bf Proof of Theorem \ref{RegSC}.}
(i)
Without loss of generality we can assume $a(z)$
has compact support because if $a(z)=0$, say for
$|z|\le 2R$, then the result is obvious.

Let $f=[S,a]u$.
We will write $C_1, C_2$, ... for constants
that may depend on $\alpha$, $\beta$, and $R$.
First of all
\[
|f(z)| \le \|a\|_{ C^\alpha}\|u\|_\infty
\int_{|t|< R} |t-z|^{\alpha-2}|d^2t|.
\]
Hence $\|f\|_\infty\le C_1 \|a\|_{ C^\alpha}\|u\|_\infty$.
Define
\[
K(z,t)=K(z)(t)=\frac{a(t)-a(z)}{(t-z)^2}.
\]
In estimating $\Delta_h f(z)$ for simplicity put $z=0$.
Introduce
\[
F(t)=\Delta_h K(0)(t)=\frac{a(t)-a(h)}{(t-h)^2}
-\frac{a(t)-a(0)}{t^2}.
\]
For $|t|<2|h|$ we use the estimate
\[
|F(t)|\le \|a\|_{ C^\alpha} (|t-h|^{\alpha-2}
+|t|^{\alpha-2}).
\]
For $|t|>2|h|$ we rewrite $F(t)$ in the form
\[
F(t)= \frac{h(2t-h)(a(t)-a(h))}{t^2(t-h)^2}
-\frac{a(h)-a(0)}{t^2},
\]
in which the first term does not exceed
$C_2 \|a\|_{ C^\alpha}|t|^{\alpha-3}|h|$.

We have
$\Delta_h f(0)=\int_\C F(t)u(t)d^2t
=J_1+J_2-(a(h)-a(0))J_3$. Here
\begin{align}
& J_1=\int_{|t|< 2|h|}F(t)u(t)\,d^2t,\qquad
J_2=\int_{|t|> 2|h|}\frac{h(2t-h)(a(t)-a(h))}{t^2(t-h)^2}
u(t)\,d^2t,\label{23}
\\
&J_3=\int_{|t|> 2|h|}\frac{u(t)}{t^2}\,d^2t.\notag
\end{align}
Then $J_1$ and $J_2$ admit the following similar estimates
\begin{align*}
&|J_1|\le \|a\|_{ C^\alpha}\left|\int_{|t|< 2|h|}
(|t-h|^{\alpha-2}+|t|^{\alpha-2})
u(t)\,d^2t\right|
\le C_3 \|a\|_{ C^\alpha}\|u\|_\infty|h|^\alpha,
\\
&|J_2|\le C_4 \|a\|_{ C^\alpha}|h|\left|\int_{|t|> 2|h|} |t|^{\alpha-3}
u(t)\,d^2t\right|
\le C_5 \|a\|_{ C^\alpha}\|u\|_\infty|h|^\alpha.
\end{align*}
Let $u\in L^\infty_R$. Then $J_3$ has the obvious estimate
\[
|J_3|\le \|u\|_\infty \int_{2|h|<|t|<R} |t|^{-2}|d^2t|
\le C_6\log|h|^{-1}\,\|u\|_\infty.
\]
Hence $\|\Delta_h f\|_\infty\le C_7 \|a\|_{ C^\alpha}
\|u\|_\infty|h|^\beta$ and $f\in C^\beta(\C)$,
which completes the proof of the first assertion in part (i).

Let $u\in C^\beta_R$. Since $S$ is bounded in $C^\beta$,
\begin{equation*}
|J_3|=\left|Su(0)
-\int_{|t|<2|h|}\frac{u(t)-u(0)}{t^2}\,d^2t\right|
\le |Su(0)|+C_6\|u\|_{C^\beta}|h|^\beta
\le C_7 ||u||_{C^\beta}.
\end{equation*}
Hence $f\in C^\alpha(\C)$, which completes the proof of (i).

(ii) Let $k=0$, $a\in C^{1,\alpha}(\C)$, $u\in C^\alpha_R$,
and $f=[S,a]u$.
We estimate the second difference $\Delta^2_hf$.
It suffices to show
$|\Delta^2_hf|\le C_8\|a\|_{ C^{1,\alpha}} \, \|u\|_{C^\alpha}|h|^{1+\alpha}$.
By Lemma \ref{Delta-Commutator}
\[
\Delta^2_hf
=\Delta_h[S,\Delta_h a]\,\delta_h u
+\Delta_h[S,a]\,\Delta_h u.
\]
Consider the first term
$A_1=\Delta_h[S,b]\,v$,
here $b=\Delta_h a$, $v=\delta_h u$.
Then by the same method we obtain
$|A_1(0)|\le C_8 \|b\|_{C^\alpha}
\|v\|_{C^\alpha}|h|^\alpha$.
Clearly $\|b\|_{C^\alpha}\le C_9 \|a\|_{ C^{1,\alpha}}|h|$
and $\|v\|_{C^\alpha}=\|u\|_{C^\alpha}$.
Hence $\|A_1\|_\infty\le C_{10} \|a\|_{ C^{1,\alpha}}\,
\|u\|_{C^\alpha}|h|^{1+\alpha}$.

The second term
$A_2=\Delta_h[S,a]\,v$ is more involved.
Here $v=\Delta_h u$, $\|v\|_\infty\le \|u\|_{C^\alpha}|h|^{\alpha}$.
Using the same notation as above, we write
$|A_2(0)|\le |J_1|+|J_2|+\|a\|_{ C^{1,\alpha}} |h J_3|$.
The terms $J_1$ and $J_3$ are handled in the same manner
as above; they admit the desired estimate. In particular,
\begin{equation*}
|J_3|
=\left|\int_{|t|>2|h|}\frac{\Delta_hu(t)}{t^2}\,d^2t\right|
=\left|\Delta_h Su(0)
-\int_{|t|<2|h|}\frac{\Delta_hu(t)-\Delta_hu(0)}{t^2}\,d^2t\right|
\le C_{11} \|u\|_{C^\alpha}|h|^\alpha.
\end{equation*}
We rewrite the remaining term $J_2=J_4+ J_5$
as a result of splitting the factor $(2t-h)$ in \eqref{23}
into the sum $2t-h=h+2(t-h)$. Then
\[
J_4=h^2
\int_{|t|> 2|h|}\frac{a(t)-a(h)}{t^2(t-h)^2}
v(t)\,d^2t,\qquad
J_5=2h
\int_{|t|> 2|h|}\frac{a(t)-a(h)}{t^2(t-h)}
v(t)\,d^2t.
\]
Since $\|v\|_\infty\le \|u\|_{C^\alpha} |h|^\alpha$,
the integral $J_4$ admits a simple estimate
\[
|J_4|\le \|a\|_{ C^{1,\alpha}} \|u\|_{C^\alpha}
|h|^{2+\alpha}\int_{|t|>2|h|}|t|^{-3}\,|d^2t|
\le C_{12}\|a\|_{ C^{1,\alpha}} \|u\|_{C^\alpha}|h|^{1+\alpha}.
\]
For the remaining term $J_5$, we use Taylor's formula
\begin{equation}
\label{21}
a(t)-a(z)=a_z(z)(t-z)+a_{\bar z}(z)\bar{(t-z)}
+O(|t-z|^{1+\alpha}).
\end{equation}
Then
$J_5=2h(a_z(h)J_6+a_{\bar z}(h)J_7+J_8)$, here
\begin{equation*}
J_6
=\int_{|t|> 2|h|}\frac{\Delta_hu(t)\,d^2t}{t^2},
\qquad
J_7
=\int_{|t|> 2|h|}\frac{\bar{t-h}}{t-h}\,
\frac{\Delta_hu(t)\,d^2t}{t^2},
\end{equation*}
and $J_8$ comes from the remainder in \eqref{21}.
The term $J_8$ has the order $|h|^{2\alpha}$,
which is even better that we need.
The term $J_6$ is the same as $J_3$ above.
Hence the desired result
for the commutator $[S,a]u$ is equivalent to the estimate
\begin{equation}
\label{15}
|J_7|\le C_{13}\|u\|_{C^\alpha}|h|^\alpha,
\end{equation}
which is independent of $a$.
Instead of dealing with $J_7$ directly, we observe
that \eqref{15} is equivalent to the desired result for
$[S,a]u$ with $a(z)=\bar z$. In this case
the commutator turns into an integral similar to the well
known Cauchy-Green operator \eqref{17}, for which the needed
result is well known (see \cite{AIM, V}). This remark completes
the proof of (ii) for $k=0$.

We now consider $k>0$. By induction we assume that the result
is already known for lower values of $k$.
Let $D_hu(z)=\frac{d}{dt}\big|_{t=0}u(z+th)$ denote
the directional derivative of $u$ in the direction $h$.
Let $a\in C^{k+1,\alpha}$ and $u\in C^{k,\alpha}_R$. Then
the result obtained for $k=0$ lets us pass to the limit
in \eqref{3} to obtain
\begin{equation*}
D_h [S,a]u=[S,D_h a]\,u +[S,a]\,D_h u.
\end{equation*}
Now by induction the result holds for all $k\ge0$.
The proof of Theorem \ref{RegSC} is complete.
$\blacksquare$

\begin{remark}
\label{RegSC-remark}
{\rm
The commutator $[S,a]$ in Theorem \ref{RegSC}(i) is in fact
a bounded operator
$L^p_R\to C^\beta(\C)$ for $p=\frac{2}{\alpha-\beta}$.
Indeed, along the lines of the above proof one can show that
$\|K(z)\|_q\le C \|a\|_{C^\alpha}$ and
$\|\Delta_h K(z)\|_q\le C \|a\|_{C^\alpha} |h|^\beta$,
here $\frac{1}{p}+\frac{1}{q}=1$.
Then the conclusion follows by H\"older inequality.
}
\end{remark}

\section{Commutators in a bounded domain}

We extend the result of the previous section to
a bounded domain.
\begin{thm}
\label{RegSD}
Let $\Omega\subset\C$ be a bounded domain of class $C^{1,\alpha}$,
$0<\alpha<1$.
\begin{itemize}
\item[(i)] If $a\in C^\alpha(\Omega)$,
then for every $0<\beta<\alpha$, the commutator $[S_\Omega,a]$
is a bounded operator $L^\infty(\Omega)\to C^\beta(\Omega)$
and $C^\beta(\Omega)\to C^\alpha(\Omega)$.
\item[(ii)]
If $\Omega$ and $a(z)$ are smooth of class $C^{k+1,\alpha}$,
$k\ge0$, then the commutator $[S_\Omega,a]$ is a bounded
operator $C^{k,\alpha}(\Omega)\to C^{k+1,\alpha}(\Omega)$.
\end{itemize}
\end{thm}

Taking into account Remark \ref{RegSC-remark},
the commutator $[S_\Omega,a]$ in (i)
is in fact a bounded operator $L^p(\Omega)\to C^\beta(\Omega)$
for $p=\frac{2}{\alpha-\beta}$.

We first recall some simple estimates.
Denote by $s=\dist(z,b\Omega)$ the distance from $z$ to $b\Omega$.
Let $r, n\ge0$ be integers.
Introduce
\[
Q_n^r u(z)=\int_{\C\setminus\Omega}
\frac{(\bar t-\bar z)^r u(t)\,dt\wedge d\bar t}{(t-z)^n}.
\]

\begin{lemma}
\label{lemma6}
Let $\Omega\subset\C$ be a bounded domain of class
$C^{k+1,\alpha}$, $k\ge0$, $0<\alpha<1$.
Let $u\in C^{k,\alpha}(\C)$.
Then there is a constant $C>0$ depending on $\Omega$, $k$,
and $\alpha$ so that for $z\in\Omega$
\begin{equation}
\label{19}
|Q_n^r u(z)| \le
\begin{cases}
C\|u\|_{C^{k,\alpha}} &\text{if}\;\; 3\le n-r< k+3,\\
C\|u\|_{C^{k,\alpha}}s^{\alpha-1} &\text{if}\;\; n-r=k+3.
\end{cases}
\end{equation}
\end{lemma}
\proof
Using induction on $k$, let $k=0$, $n-r=3$.
Then
\[
Q_n^ru(z)=\int_{\C\setminus\Omega}
\frac{(\bar t-\bar z)^r(u(t)-u(z))\,dt\wedge d\bar t}
{(t-z)^n}
+u(z)Q_n^r(1)(z),
\]
in which the first term clearly admits the estimate
$O(s^{\alpha-1})$ by integrating the modulus of the integrand.
The second term will be automatically
considered simultaneously with the general case.

Now let $k\ge1$ and assume the estimate \eqref{19} for lower values
of $k$. We also allow $k=0$, $u\equiv1$.
Let $b\Omega$ be a level set of a function of class
$C^{k+1,\alpha}$. Then on $b\Omega$ we have $d\bar t=\phi(t)dt$,
here $\phi\in C^{k,\alpha}$. We assume $\phi$ extends
to the whole plane and has compact support.
Introduce
\begin{equation}
\label{31}
K_m^r u(z)=\int_{b\Omega}
\frac{(\bar t-\bar z)^r u(t)\,dt}{(t-z)^m}.
\end{equation}
By Stokes' formula
\begin{align*}
K_{n-1}^r(u\phi)=(n-1)Q_n^r u+Q_{n-1}^r u_z.
\end{align*}
The term $Q_{n-1}^r u_z$ satisfies \eqref{19}
by induction.
Hence it suffices to show that $K_{n-1}^r(u)$
satisfies \eqref{19} for $u\in C^{k,\alpha}(\C)$,
$k\ge0$.

Integrating by parts for $m>1$ yields
\begin{equation*}
(1-m)K_m^r u=rK_{m-1}^{r-1}(u\phi)+
K_{m-1}^r(u_z+u_{\bar z}\phi).
\end{equation*}
Starting with $m=n-1$, we successively integrate
by parts all resulting terms while still possible.
If $n-r<k+3$, then all final terms will have the
form $K_m^p(\psi)$ with $m<p+2$, $\psi\in C^\alpha$.
They are clearly bounded.
If $n-r=k+3$, then all final terms will have the
form $K_{p+2}^p(\psi)$ with $\psi\in C^\alpha$.
We have
\[
K_{p+2}^p(\psi)(z)
=\int_{b\Omega}
\frac{(\bar t-\bar z)^p(\psi(t)-\psi(z))\,dt\wedge d\bar t}
{(t-z)^{p+2}}
+\psi(z)K_{p+2}^p(1)(z),
\]
The first term clearly has the estimate $O(s^{\alpha-1})$.
For the second one, we again reduce the exponent $p$ by integration
by parts and eventually obtain the same estimate $O(s^{\alpha-1})$.
$\blacksquare$
\medskip

\noindent
{\bf Proof of Theorem \ref{RegSD}.}
(i) The assertion about the map
$L^\infty(\Omega)\to C^\beta(\Omega)$
follows immediately from the corresponding assertion
of Theorem \ref{RegSC} because $a\in C^\alpha(\Omega)$
can be extended to a function of $C^\alpha(\C)$, and
for $R$ large enough, $L^\infty(\Omega)\subset L^\infty_R$
by trivial extension. So we focus on $[S_\Omega,a]$
on the space $C^\beta(\Omega)$.

We again consider the extension of $a\in C^\alpha(\Omega)$
to the whole plane (which we still denote
by the same letter $a$) that has compact support,
smooth in $\C\setminus\bar\Omega$, and whose first derivatives
admit the estimate $O(s^{\alpha-1})$.
Similarly, we extend the restriction
$a|_{b\Omega}$ inside $\Omega$ so that the extension $\tilde a$ is
smooth in $\Omega$ with first derivatives of the magnitude
$O(s^{\alpha-1})$. Finally, we extend the given function
$u\in C^\beta(\Omega)$ to the whole plane
so that the extension has compact support and belongs to
$C^\beta(\C)$.
Then for $z\in\Omega$,
\begin{align*}
&[S_\Omega,a]u(z)=[S_\C,a]u(z)-v_1(z)-v_2(z),\\
&v_1(z)=\int_{\C\setminus\Omega} \frac{a(t)-\tilde a(z)}{(t-z)^2}\,
u(t)\,d^2t, \quad
v_2(z)=\int_{\C\setminus\Omega} \frac{\tilde a(z)-a(z)}{(t-z)^2}\,
u(t)\,d^2t.
\end{align*}
By Theorem \ref{RegSC}, the first term $[S_\C,a]u$ has the desired
properties. Differentiating $v_1$ yields
\begin{equation}
\label{1}
(v_1)_z(z)=-\tilde a_z(z)
\int_{\C\setminus\Omega} \frac{u(t)\,d^2t}{(t-z)^2}
+2\int_{\C\setminus\Omega} \frac{a(t)-\tilde a(z)}{(t-z)^3}\,
u(t)\,d^2t.
\end{equation}
For $u\in C^\beta$
the first integral in (\ref{1}) is bounded.
Since $\tilde a_z=O(s^{\alpha-1})$, the first term
in (\ref{1}) is $O(s^{\alpha-1})$.

The second integral in (\ref{1}) has the estimate
$O(\int_{|t|>s}|t|^{\alpha-3}|d^2t|)=O(s^{\alpha-1})$.
The $\bar z$-derivative of $v_1$ is estimated similarly
but slightly simpler.
Then by Hardy-Littlewood lemma, $v_1\in C^\alpha(\Omega)$.

We now consider $v=v_2=bw$. Here $b=\tilde a-a\in C^\alpha(\Omega)$,
$b|_{b\Omega}=0$, $b=O(s^\alpha)$; $w=S_{\C\setminus\Omega}u$.
Let $z,z'\in \Omega$; without loss of generality $z$ is closer
to $b\Omega$ than $z'$. We estimate $\Delta v=v(z)-v(z')$ in terms
of $h=|z-z'|$.
We have
\begin{equation}
\label{2}
\Delta v=\Delta b\, w(z)+b(z')\Delta w.
\end{equation}
Since $u\in C^\beta(\C)$,
we have $w(z)=O(1)$, $w_z(z)=O(s^{-1})$,
$\Delta w=O(hs^{-1})$.
Let $h\le s$. Then $\Delta b=O(h^\alpha)$,
$b(z')=O(s^\alpha)$. Plugging these estimates in
(\ref{2}) yields
\[
\Delta v=O(h^\alpha +s^\alpha hs^{-1})=O(h^\alpha).
\]
If $h\ge s$, then $v(z)=O(s^\alpha)$,
$v(z')=O(h^\alpha)$, and again $\Delta v=O(h^\alpha)$.
Hence $[S_\Omega,a]u\in C^\alpha(\Omega)$ as desired.

(ii) Let $a\in C^{k+1,\alpha}(\Omega)$.
We again assume that $a(z)$ and $u(z)$ are extended
to the whole plane.
(We do not need the other extension $\tilde a$.)
We represent
\[
[S_\Omega,a]u=[S_\C,a]u-v,\qquad
v=[S_{\C\setminus\Omega},a]u,\qquad
v(z)=\int_{\C\setminus\Omega} \frac{a(t)-a(z)}{(t-z)^2}\,
u(t)\,d^2t.
\]
By Theorem \ref{RegSC}, $[S_\C,a]u\in C^{k+1,\alpha}(\Omega)$.
To see that $v\in C^{k+1,\alpha}(\Omega)$, we first differentiate
it $(k+1)$ times. One term will have the form
\[
-D^{k+1}a \, S_{\C\setminus\Omega}u,
\]
here $D^{k+1}a$ denotes any derivative of order $(k+1)$.
This term is clearly in $C^\alpha(\Omega)$.
To show that the other terms are in $C^\alpha(\Omega)$,
we show that the first derivatives of these terms have
the estimate $O(s^{\alpha-1})$.
Then by Hardy-Littlewood lemma we will obtain
$v\in C^{k+1,\alpha}(\Omega)$.
By differentiating one more time, we obtain the following terms.
There will be one term of the form
\[
J_1(z)=\int_{\C\setminus\Omega}\frac{a(t)-a(z)}{(t-z)^{k+4}}
u(t)\,d^2t,
\]
while all other terms will be constant multiples of
integrals of the form
\[
J_2(z)=D^p a(z)
\int_{\C\setminus\Omega}\frac{u(t)\,d^2t}{(t-z)^{q+2}}.
\]
Here $p+q=k+2$, $p\ge1$, hence $q+2\le k+3$. Then the terms
of the form $J_2$ are all bounded by Lemma \ref{lemma6}.
For the term $J_1$, we use Taylor's formula
\begin{equation}
\label{22}
a(t)-a(z)=\sum_{1\le p+q\le k+1}a_{pq}(z)(t-z)^p(\bar t-\bar z)^q
+O(|t-z|^{k+1+\alpha}).
\end{equation}
The term corresponding to the remainder in \eqref{22}
is estimated directly; it has the order $O(s^{\alpha-1})$.
Now by Lemma \ref{lemma6} the estimate $J_1(z)=O(s^{\alpha-1})$
follows. Theorem is proved.
$\blacksquare$
\medskip

Finally we include a simple result that applies to
$[S_\Omega,a]$ above.
\begin{prop}
\label{Increase-p}
Let $\Omega\subset\C$ be a bounded domain.
Let $k:\Omega\times\Omega\to\C$ satisfy
$|k(z,t)|\le |z-t|^{\alpha-2}$, $0<\alpha<1$.
Let $Ru(z)=\int_\Omega k(z,t)u(t)\,d^2t$.
Let $p\ge1$.
Let $p\le r<\frac{2p}{2-\alpha p}$
(if $2-\alpha p<0$, then $p\le r\le\infty$).
Then $R:L^p(\Omega)\to L^r(\Omega)$ is a bounded operator,
which in particular holds if $r=p+\frac{\alpha}{2}$.
\end{prop}
\proof
For completeness we include a proof.
Since $\Omega$ is bounded, for $\lambda>-2$
there is a constant $C(\lambda)>0$ such that
$\int_\Omega|z-t|^\lambda|d^2t|\le C(\lambda)$.

Let $0\le c\le 1$, $\frac{1}{p}+\frac{1}{q}=1$,
$\frac{1}{r}+\frac{1}{s}+\frac{1}{q}=1$,
hence $s=\frac{rp}{r-p}$. By H\"older inequality,
\begin{align*}
|Ru(z)|
&=\left|\int_\Omega k(z,t)u(t)\,d^2t\right|\le
\int_\Omega |u|^{\frac{p}{r}} |k|^c |u|^{1-\frac{p}{r}}
|k|^{1-c}\, |d^2t| \\
&\le
\left(\int_\Omega |u|^p |k|^{cr}\,|d^2t|\right)^{1/r}
\left(\int_\Omega |u|^p\,|d^2t|\right)^{1/s}
\left(\int_\Omega |k|^{(1-c)q}\,|d^2t|\right)^{1/q},\\
\|Ru\|^r_r
&=\int_\Omega |Ru(z)|^r\,|d^2z|
\le C((\alpha-2)cr) C((\alpha-2)(1-c)q)^{r/q} \|u\|^r_p,
\end{align*}
provided that
$(\alpha-2)cr>-2$ and $(\alpha-2)(1-c)q>-2$.
These conditions yield the desired bounds for $r$.
$\blacksquare$

\section{The operator $S_\Omega\bar S_\Omega$}

The standard approach to the Beltrami type equations
(see \cite{AIM, V}) involves the Cauchy-Green operator
\begin{equation}
\label{17}
T_\Omega u(z)=\int_\Omega \frac{u(t)\,d^2t}{t-z}
\end{equation}
for a domain $\Omega\subset\C$. We also consider
its modification suitable for solving the Dirichlet
problem in the unit disc $\D=\{z\in\C:|z|<1\}$.
\[
T_1 u(z)=T_\D u(z)-\bar{T_\D u(\bar z^{-1})}
=\int_\D \frac{u(t)\,d^2t}{t-z}
+\int_\D \frac{z \bar {u(t)}\,d^2t}{1-z\bar t}.
\]
Both operators $T_\Omega$ and $T_1$ solve the
$\bar\partial$-problem, and $T_1$ in addition satisfies
the boundary condition $\Re T_1u|_{b\D}=0$.
Consider the operators $S_\Omega$ and $S_1$,
the $\partial$-derivatives of $T_\Omega$ and $T_1$.
Then
\[
S_1u(z)=\pv\int_\D \frac{u(t)\,d^2t}{(t-z)^2}
+\int_\D \frac{\bar {u(t)}\,d^2t}{(1-z\bar t)^2}.
\]
In other words
\[
S_1u=S_\D u-B\bar u,\qquad
B v(z)=-\int_\D \frac{v(t)\,d^2t}{(1-z\bar t)^2}.
\]
Here $B$ is the Bergmann projection in $\D$
(in this notation $d^2t<0$).
It is well known (see \cite{AIM, V}) that $S_\C$ and $S_1$
are isometries of $L^2(\C)$ and $L^2(\D)$ respectively, that is,
$S_\C \bar{S_\C}=I$ and $S_1 \bar{S_1}=I$, here $I=\id$
is the identity operator.
(Note $S_\C^*=\bar{S_\C}$ and $S_1^*=\bar{S_1}$.)
We make the following observation.

\begin{lemma}
\label{lemma4}
$S_\D \bar{S_\D}=I-B$ and $BS_\D+S_\D\bar B=0$.
\end{lemma}
\proof
Introduce the conjugation operator
$\iota u=\bar u$. Then $\iota^2=I$.
For every operator $P$ we have by definition
$\bar P=\iota P\iota$ or $\iota P=\bar P\iota$,
in particular $\bar\iota=\iota$.
For simplicity put $S=S_\D$. Then we write
$S_1=S-B\iota$, $\bar S_1=\bar S-\bar B\iota$.
Using $B^2=B$ we obtain
\[
\id=S_1 \bar{S_1}
=(S-B\iota)(\bar S-\bar B\iota)
=S\bar S-B\iota\bar S-S\bar B\iota+B\iota\bar B\iota
=(S\bar S+B)-(BS+S\bar B)\iota.
\]
Now by separating linear and anti-linear terms, we
obtain the desired relations.
$\blacksquare$
\medskip

We now consider the operator
$B_\Omega=I-S_\Omega\bar S_\Omega$ for an arbitrary
smooth domain $\Omega\subset\C$.
We will see in the next section that $B_\Omega$ is
related to the Bergmann projection for $\Omega$.
Here we only care to what extent $B^2=B$ holds for $B_\Omega$.
Invoke the Cauchy type integral
\begin{equation}
\label{41}
K_\Omega u(z)=\frac{1}{2\pi i}\int_{b\Omega}\frac{u(t)\,dt}{t-z},
\qquad
z\in \Omega.
\end{equation}
For $z\in b\Omega$ we interpret $K_\Omega u(z)$ as a boundary
value of the function $K_\Omega u$ in $\Omega$.
With some abuse of notation we write
$\partial u(z)= \partial_z u(z)$ and
$\bar\partial u(z)= \partial_{\bar z} u(z)$.
We recall the Cauchy-Green-Pompeiu formula
\[
K_\Omega+T_\Omega\bar\partial=I.
\]
The following result is similar to one by
Kerzman and Stein \cite{KS} who discovered that
$K_\Omega-K_\Omega^*$ is a smoothing operator.

\begin{thm}
\label{Kerzman-Stein}
Let $\Omega\subset\C$ be a bounded domain of class
$C^k$, here $k\ge 1$ may be fractional.
Then $P=K_\Omega+\bar K_\Omega-I$
is a bounded operator
$L^1(b\Omega)\to C^{k-2}(b\Omega)$ (if $k\ge2$) and
$C^1(b\Omega)\to C^{k-1}(b\Omega)$.
\end{thm}
\proof
Observe for $z\in\Omega$, $t\in b\Omega$
\[
\frac{dt}{t-z}-\frac{d\bar t}{\bar t-\bar z}
=\frac{\bar t-\bar z}{t-z}\,\,d_t \(\frac{t-z}{\bar t-\bar z}\)
=d_t\log\(\frac{t-z}{\bar t-\bar z}\)
=2i\,d_t\arg(t-z).
\]
Then for $z\in\Omega$
\[
(K_\Omega+\bar K_\Omega)u(z)
=\frac{1}{2\pi i}\int_{b\Omega}u(t)
\(\frac{dt}{t-z}-\frac{d\bar t}{\bar t-\bar z}\)
=\frac{1}{\pi}\int_{b\Omega}u(t) \,d_t\arg(t-z).
\]
If $z\in b\Omega$, then passing to the limit yields
\[
(K_\Omega+\bar K_\Omega)u(z)
=u(z)+\frac{1}{\pi}\,\pv\int_{b\Omega}u(t) \,d_t\arg(t-z).
\]
The only reason for principal value in this integral
is the jump of $-\pi$ of $\arg(t-z)$ at $t=z$.
Otherwise the integral has a smooth kernel. Indeed,
suppose an arc of $b\Omega$ has a parametric equation
$\gamma(\tau)=\tau+i\phi(\tau)$ with $\phi\in C^k$.
Then for $t=\gamma(\tau)$, $z=\gamma(\tau_0)$, we have
\[
d_t\arg(t-z)
=d\arctan \(\frac{\phi(\tau)-\phi(\tau_0)}{\tau-\tau_0}\)
\in C^{k-2}.
\]
Hence if $u\in L^1(b\Omega)$, then $Pu\in C^{k-2}(b\Omega)$.
If $u\in C^1(b\Omega)$, then by integrating by parts
$Pu\in C^{k-1}(b\Omega)$.
$\blacksquare$
\medskip

\begin{cor}
\label{B-square}
Let $\Omega\subset\C$ be a bounded domain of class
$C^k$, here $k>2$ is fractional. Then
$B_\Omega^2-B_\Omega$ is a bounded operator
$L^p(\Omega)\to C^{k-3}(\bar\Omega)$ (if $p>1$, $k>3$) and
$C^\alpha(\Omega)\to C^{k-2}(\bar\Omega)$
(if $0<\alpha<1$, $k>2$).
\end{cor}
\proof
For simplicity of notation we omit the subscripts $\Omega$.
We have
\begin{align}
& KP=K(K+\bar K-I)=K+K\bar K-K=K\bar K,\notag\\
& B=I-S\bar S=I-\partial (T\bar\partial)\bar T
=I-\partial(I-K)\bar T=\partial K\bar T,\label{44}\\
& B^2=\partial K(\bar T \partial) K\bar T
=\partial K(I-\bar K) K\bar T
=B-\partial KPK\bar T.\notag
\end{align}
Let $u\in L^p(\Omega)$, $p>1$.
Then $\bar T u\in W^{1,p}(\Omega)$,
the trace $\bar T u|_{b\Omega}\in L^p(b\Omega)$,
and $K\bar Tu\in L^p(b\Omega)$.
By Theorem \ref{Kerzman-Stein},
$PK\bar Tu\in C^{k-2}(b\Omega)$.
Since $k$ is fractional,
we have $KPK\bar Tu\in C^{k-2}(\Omega)$,
and $(B^2-B)u=-\partial KPK\bar Tu\in C^{k-3}(\Omega)$,
as desired.

Let $u\in C^\alpha(\Omega)$.
Then $\bar T u\in C^{1+\alpha}(\Omega)$ and
$K\bar Tu\in C^{1+\alpha}(\Omega)$.
By Theorem \ref{Kerzman-Stein},
$PK\bar Tu\in C^{k-1}(b\Omega)$.
Since $k$ is fractional, we have
$KPK\bar Tu\in C^{k-1}(\Omega)$,
and $(B^2-B)u\in C^{k-2}(\Omega)$, as desired.
$\blacksquare$

\section{The Bergman projection}

The main result of this section is an asymptotic formula for
the Bergman projection ${\cal B}={\cal B}_\Omega$ for
a smooth bounded simply connected domain $\Omega\subset \C$.
We will not need it in the rest of the paper.

Let $H=H(\Omega)$ be the Bergman space of all holomorphic
functions of class $L^2(\Omega)$. The Bergman projection
${\cal B}:L^2(\Omega)\to H$ is the orthogonal projection
onto the subspace $H\subset L^2(\Omega)$.

\begin{thm}
\label{Bergman}
Let $\Omega\subset\C$ be a bounded domain of class
$C^k$, here $k>3$ is fractional.
Let $B_n=I-(S\bar S)^n$, $S=S_\Omega$, $n\ge 1$.
Then
\begin{itemize}
\item[(i)]
For all $n\ge1$, the difference ${\cal B}-B_n$
is a bounded operator $L^2(\Omega)\to C^{k-3}(\Omega)$.
\item[(ii)] If $\Omega$ is simply connected, then
${\cal B}=\lim_{n\to\infty}B_n$.
\end{itemize}
\end{thm}

We can compare this result to the one by Kerzman and Stein \cite{KS}.
Let ${\cal S}:L^2(b\Omega)\to H^2(b\Omega)$ be the Szeg\"o
orthogonal projection, and let $K=K_\Omega$ be the Cauchy transform.
Kerzman and Stein \cite{KS} proved that ${\cal S}-K$ is a compact
smoothing operator and ${\cal S}=K(I-A)^{-1}$, here $A=K-K^*$
is a compact smoothing operator. If $\Omega$ is sufficiently close
to the disc $\D$, then $\|A\|<1$, and the inverse has an explicit
formula $(I-A)^{-1}=\sum_{n=0}^\infty A^n$. In contrast,
our formula ${\cal B}=\lim_{n\to\infty}B_n$ holds for every simply
connected smooth domain.

We also compare Theorem \ref{Bergman}(i) with the formula
${\cal B}=\partial E \bar T$ (Bell \cite{Bell}, page 70).
Here $E$ denotes the harmonic extension from $b\Omega$
to $\Omega$. For a general domain, $E$ is not explicit.
If we replace $E$ by $K$, then by \eqref{44} we obtain
the explicit operator $B=B_1$, which by Theorem \ref{Bergman}(i) approximates the Bergman projection $\cal B$.

\begin{lemma}
\label{lemma8}
The subspace $H\subset L^2(\Omega)$ is invariant for $S\bar S$,
and $S\bar S|_H: H\to C^{k-3}(\bar\Omega)$ is bounded.
\end{lemma}
\proof
By \eqref{44},
the subspace $H\subset L^2(\Omega)$ is invariant for $B=I-S\bar S$,
hence for $S\bar S$.
Let $H_0\subset H$ consist of such $u\in H$ that for every
closed path $\gamma\subset\Omega$, we have $\int_\gamma u(z)\,dz=0$.
For $u\in H_0$ define $Ju(z)=\int_{z_0}^z u(t)\,dt$ along a path
in $\Omega$. Then for $u\in H_0$
\[
S\bar Su=\partial(T\bar\partial)(\bar T\partial)Ju
=\partial(I-K)(I-\bar K)Ju
=\partial(I-K-\bar K+K\bar K)Ju
=\partial K\bar KJu
=\partial KPJu.
\]
Since $u\in L^2(\Omega)$, we have $Ju\in W^{1,2}(\Omega)$.
Then the trace $Ju|_{b\Omega}\in L^2(b\Omega)$.
Then by Theorem \ref{Kerzman-Stein}, we have
$PJu\in C^{k-2}(b\Omega)$, $KPJu\in C^{k-2}(\bar\Omega)$,
and finally $S\bar Su=\partial KPJu\in C^{k-3}(\bar\Omega)$.

If $\Omega$ is simply connected, then $H_0=H$, and the
proof is complete. Otherwise, $H_0$ has a finite dimensional
(not necessarily orthogonal) complement $H_1$ in $H$
of the form, say
\[
H_1=\left\{\sum c_j(z-z_j)^{-1}: c_j\in\C \right\}.
\]
Here the points $z_j\in\C\setminus\Omega$ are fixed ---
one in each bounded component of $\C\setminus\Omega$.
Since $H_1$ consists of smooth functions,
the operator $S\bar S|_{H_1}:H_1\to C^{k-3}(\bar\Omega)$
is bounded, hence the desired conclusion.
$\blacksquare$

\begin{lemma}
\label{lemma9}
If $\Omega$ is simply connected, then
$\|S\bar S\|_H<1$.
\end{lemma}
\proof
Note $\bar S=S^*$. Since $S\bar S|_H$ is self-adjoint,
compact, and $S\bar S\ge0$, it suffices to show
that $S\bar S|_H$ does not have the eigenvalue 1.

Suppose there is $u\in H$ such that $S\bar Su=u$.
Since $\|S\|_2\le 1$, we have
$\|S\bar Su\|_2 \le \|\bar Su\|_2\le\|u\|_2$.
Since $S\bar Su=u$, we have in particular,
$\|\bar Su\|_2=\|u\|_2$.
On the other hand, $S_\C$ is an isometry of $L^2(\C)$.
Hence $\bar Su(z)=0$ for $z\notin \Omega$,
that is, $\bar\partial\bar Tu=0$ in $\C\setminus\Omega$.

Since $\bar Tu$ is antiholomorphic on a connected set
$\C\setminus\Omega$ and $\bar\partial\bar Tu=0$,
the function $\bar Tu=\const$ in $\C\setminus\Omega$.
In fact $\bar Tu|_{\C\setminus\Omega}=0$ because
it vanishes at infinity.
By Lemma \ref{lemma8}, $u\in C^{k-3}(\bar\Omega)$, hence
$\bar Tu$ is continuous on $\C$.

Since $u$ is holomorphic,
$\bar\partial\partial\bar Tu=\bar\partial u=0$,
that is, $\bar Tu$ is harmonic in $\Omega$.
Since $\bar Tu|_{b\Omega}=0$, we have $\bar Tu=0$ in $\Omega$.
Hence $u=\partial \bar Tu=0$, and the proof is complete.
$\blacksquare$
\bigskip

\noindent
{\bf Proof of Theorem \ref{Bergman}.}
By \eqref{44}, we have $B(L^2(\Omega))\subset H$.
Since $B$ is self-adjoint, $B(H^\bot)=0$.
Indeed, for every $u\in H^\bot$ and $v\in L^2(\Omega)$,
we have $(Bu, v)=(u, Bv)=0$.
Hence, $S\bar S H\subset H$ and $S\bar S|_{H^\bot}=I$.

We now compare $B_n$ with $\cal B$ on $H$ and $H^\bot$.
On $H^\bot$ we have ${\cal B}(H^\bot)=B_n (H^\bot)=0$.
On $H$ we have
$({\cal B}-B_n)|_H=(S\bar S)^n|_H:H\to C^{k-3}(\bar \Omega)$,
which proves (i).

By Lemma \ref{lemma9}, $(S\bar S)^n|_H\to 0$ as $n\to\infty$,
hence the conclusion (ii).

Theorem is proved.
$\blacksquare$
\medskip

We realize that if $\Omega$ is not simply connected,
then Lemma \ref{lemma9} and Theorem \ref{Bergman}(ii) fail
as the following simple example shows.

\begin{example}
\label{example1}
{\rm
Let $0<r<1$ and let
$\Omega=\{z: r<|z|<1\}$.
Let $u(z)=1/z$. Then one can find $\bar Tu(z)=2\log|z|$.
(It is independent of $r$.)
Then
$S\bar Su(z)
=\partial T\bar\partial\bar Tu(z)
=\partial T\bar\partial (2\log|z|)
=\partial T(1/\bar z)
=\partial (2\log|z|)
=1/z
=u(z)$.
Then $B_n u=0$, but ${\cal B}u=u$,
so Lemma \ref{lemma9} and Theorem \ref{Bergman}(ii) fail.
}
\end{example}

\section{Integral equations with operator $S$}

We now consider the integral equation
\begin{equation}
\label{9}
u=S(A\bar u)+b.
\end{equation}
Here $u$ and $b$ are $m$-vector functions and $A$ is a
$m\times m$ matrix function in a smooth bounded domain
$\Omega\subset\C$; $m\ge1$, $S=S_\Omega$. In the future,
with some abuse of notation,
we omit the parentheses in (\ref{9})
and similar equations, interpreting $A$ as the operator
of multiplication by $A$.
We impose the condition $\|A\|_\infty<1$.
Here $\|A\|_\infty$ denotes the
maximum of the Euclidean operator norm of $A(z)$
over all $z\in\bar\Omega$.

\begin{prop}
\label{IntEquation-SD}
Let $\Omega\subset\C$ be a bounded domain of class $C^\infty$.
Let $A,b\in C^{k,\alpha}(\Omega)$, $0<\alpha<1$, $k\ge0$,
$\|A\|_\infty<1$.
Then the equation (\ref{9}) has a unique solution
$u\in L^2(\Omega)$. This solution
$u\in C^{k,\alpha}(\Omega)$, and for fixed $A$ the operator
$b\mapsto u$ is bounded in $C^{k,\alpha}(\Omega)$
\end{prop}
The proof below goes through if $\Omega$ has finite smoothness
of class $C^{3,\beta}$ ($0<\beta<1$) if $k=0$ and $C^{k+2,\alpha}$
if $k\ge1$.
\medskip

\proof
The existence of a unique solution
$u\in L^2(\Omega)$ is standard
(see \cite{AIM, V}). It follows because
$\|S\circ A\|_2<1$ as an operator in $L^2(\Omega)$.

Iterating (\ref{9}) yields
\[
u=SA\bar S\bar A u+b_1,\qquad
b_1=SA\bar b+b\in C^{k,\alpha}(\Omega).
\]
Interchanging $\bar S$ and $A$ yields
\[
u=S\bar S A\bar Au+b_2,\qquad
b_2=S[A,\bar S]\bar Au+b_1.
\]
We include the term $S[A,\bar S]\bar Au$ in $b_2$
because by the results of Section 3 the commutator
is ``better'' than $u$.
Recall $S\bar S=I-B$, $B=B_\Omega$. Since
$[PQ,R]=[P,R]Q+P[Q,R]$, both Theorem \ref{RegSD}
and Proposition \ref{Increase-p} apply to $B$.
We put $v=A\bar Au$. Then
\[
v=A\bar A(v-Bv)+b_3, \qquad
b_3=A\bar A b_2.
\]
Since $|A|<1$, we have $(I-A\bar A)^{-1}\in C^{k,\alpha}$, and
\[
v=-A_0 Bv+b_4,\qquad
A_0=(I-A\bar A)^{-1}A\bar A,\quad
b_4=(I-A\bar A)^{-1}b_3.
\]
Applying $B$ and interchanging $B$ and $A_0$ yields
\begin{align*}
& Bv=-B A_0 Bv+Bb_4=-([B,A_0]+A_0B)Bv+Bb_4=-A_0Bv+b_5,\\
& b_5=-[B,A_0]Bv + A_0(B-B^2)v + Bb_4.
\end{align*}
Note that by Corollary \ref{B-square} the term $(B-B^2)v$
is $C^\infty$. Also note $(I+A_0)^{-1}=I-A\bar A$ and
$A_0(I-A\bar A)=A\bar A$. Then
\begin{equation*}
\label{10}
Bv=(I-A\bar A)b_5,\qquad
v=-A\bar Ab_5+b_4,\qquad
u=S\bar Sv+b_2.
\end{equation*}
As a result, the initial equation implies
\begin{equation}
\label{25}
u=Mu+Nb,
\end{equation}
where $M$ is a smoothing operator with properties
described in Theorem \ref{RegSD} and Proposition \ref{Increase-p},
and $N$ is a bounded operator in $C^{k,\alpha}$.

We now use \eqref{25} for {\it bootstrapping}, successively
improving the regularity of the solution.
Since $u\in L^p$, starting from $p=2$,
by Proposition \ref{Increase-p},
$Mu \in L^{p+\frac{\alpha}{2}}$,
hence by \eqref{25} $u\in L^{p+\frac{\alpha}{2}}$.
We repeat this argument finitely many times
till we get $u\in L^{r}$, $r>\frac{2}{\alpha}$.
Repeating it one more time, by Proposition \ref{Increase-p}
we get $u\in L^\infty$.
We now repeat it again finitely many times using
Theorem \ref{RegSD} and get $u\in C^{k,\alpha}$ as desired.
Note that the number of times we iterate \eqref{25}
depends only on $k$ and $\alpha$.
$\blacksquare$
\medskip

We now consider a similar integral equation
in the unit disc, namely
\begin{equation}
\label{11}
u=S_1(A\bar u)+b.
\end{equation}

\begin{prop}
\label{IntEquation-S1}
For the equation (\ref{11}), Proposition \ref{IntEquation-SD} holds.
\end{prop}
\proof
The argument of the proof is similar to that for
Proposition \ref{IntEquation-SD}. The difference is that
the results of Section 3 do not directly apply to $[S_1,A]$
because $S_1$ is not complex linear. Nevertheless we
reduce the result to (the proof of) Proposition \ref{IntEquation-SD}.
By the definition of $S_1$
\[
u=SA\bar u-B\bar Au+b,
\]
here $S=S_\D$, $B=B_\D$. By Lemma \ref{lemma4}
\[
Bu=-S\bar BA\bar u-B\bar Au+Bb.
\]
Multiplying by $\bar A$ and interchanging $\bar A$ and $B$
yields
\[
B\bar Au=-\bar AS\bar BA\bar u-\bar AB\bar Au+b_1,\qquad
b_1=-[\bar A,B]u+\bar ABb.
\]
Then $v=B\bar Au$ satisfies the equation
$v=-\bar AS\bar v-\bar Av+b_1$, which in turn simplifies to
\begin{equation}
\label{12}
v=A_1 S\bar v+b_2,\qquad
A_1=-(I+\bar A)^{-1}\bar A,\quad
b_2=(I+\bar A)^{-1}b_1.
\end{equation}
This equation looks similar to (\ref{9}), however
$\|A_1\|_\infty<1$
need not hold, so the equation requires a little more care.
Following
the beginning of the proof of Proposition \ref{IntEquation-SD},
iterating (\ref{12}) yields
\[
v=A_1 S\bar A_1 \bar S v+b_3,\qquad
b_3=A_1S\bar b_2 +b_2.
\]
Interchanging $S$ and $\bar A_1$ yields
\[
v=A_1\bar A_1 S\bar Sv+b_4,\qquad
b_4=A_1[S,\bar A_1]\bar Sv+b_3.
\]
Since $v=B\bar Au$ and $B^2=B$, we have $S\bar Sv=v-Bv=0$.
Hence $v=b_4$, and the original equation takes the form
\[
u=SA\bar u+b_5,\qquad
b_5=-b_4+b,
\]
which is the subject of Proposition \ref{IntEquation-SD}.
By bootstrapping we obtain $u\in C^{k,\alpha}(\D)$.
$\blacksquare$

\section{Dirichlet problem}

We consider the Dirichlet problem for an elliptic equation
\begin{equation}
\label{18}
f_{\bar z}=a(z)f_z + b(z)\bar f_{\bar z}+c(z).
\end{equation}
In the scalar case the ellipticity means that either $|a|+|b|<1$
or $||a|-|b||>1$; the two cases are related by the interchange
$f\leftrightarrow \bar f$. We restrict to the former case.
Our main result is the following.

\begin{thm}
\label{Dirichlet-scalar}
Let $\Omega\subset\C$ be a simply connected domain
of class $C^{k+1,\alpha}$, $k\ge0$, $0<\alpha<1$.
Let $a, b, c\in C^{k,\alpha}(\Omega)$,
$f_0\in C^{k+1,\alpha}(b\Omega,\R)$,
$|a|+|b|\le a_0<1$, for some constant $a_0$.
Then the scalar equation (\ref{18}) with boundary condition
$\Re f|_{b\Omega}=f_0$ has a unique solution
in the Sobolev class $W^{1,2}(\Omega)$.
This solution $f\in C^{k+1, \alpha}(\Omega)$, and for fixed
$a$ and $b$ the map $(c,f_0)\mapsto f$ is a bounded operator
$C^{k,\alpha}(\Omega)\times C^{k+1,\alpha}(b\Omega,\R)\to
C^{k+1,\alpha}(\Omega)$.
\end{thm}

For simplicity we assume that $\Omega$ is simply connected because our method involves reduction to the unit disc.
Thus we begin the proof with several reductions.
\begin{lemma}
\label{lemma1}
It suffices to prove Theorem \ref{Dirichlet-scalar}
for $f_0=0$ and $\Omega=\D$, the unit disc.
\end{lemma}
\proof
To reduce to $f_0=0$, we fix $f_1\in C^{k+1,\alpha}(\Omega)$
satisfying $\Re f_1|_{b\Omega}=f_0$. Then for the new unknown
$\tilde f=f-f_1$, the equation will have a form similar
to the original one, and the boundary condition will
turn into $\Re \tilde f|_{b\Omega}=0$.

To reduce to $\Omega=\D$ we can introduce a new independent
variable $\zeta=\psi(z)$, so that $\psi:\bar\Omega\to\bar\D$
is a $C^{k+1,\alpha}$ diffeomorphism with positive Jacobian.
The equation will preserve its form and boundary conditions.
Moreover, if $a=0$ or $b=0$, then by choosing a conformal
map $\psi$ this condition can be preserved also.
$\blacksquare$

\begin{lemma}
\label{lemma3}
It suffices to prove Theorem \ref{Dirichlet-scalar}
for $a=0$.
\end{lemma}
\proof
We change the independent variable by a Beltrami
homeomorphism $\psi:\D\to\D$ of the equation
\begin{equation}
\label{5}
\psi_{\bar z}=\mu(z)\psi_z.
\end{equation}
The Beltrami coefficient $\mu$ will be determined later.
The equation (\ref{18}) will take the form
\begin{equation}
\label{7}
g_{\bar\zeta}=\tilde a g_\zeta
+ \tilde b\bar g_{\bar\zeta}
+\tilde c,
\end{equation}
here $g=f\circ \psi^{-1}$. We write $\zeta=\psi(z)$.
By straightforward calculations we now find the new coefficients.
We have
\[
g_{\bar z}=g_\zeta \psi_{\bar z}+g_{\bar\zeta}\bar\psi_{\bar z}
=a(g_\zeta \psi_z+g_{\bar\zeta}\bar\psi_z)
+b(\bar g_{\bar\zeta} \bar\psi_{\bar z}
+\bar g_{\zeta}\psi_{\bar z})+c.
\]
By (\ref{5}) we obtain
\begin{equation}
\label{6}
\bar\psi_{\bar z}(1-a\bar\mu)g_{\bar\zeta}
-\psi_z b\mu \bar g_\zeta
=\psi_z(a-\mu) g_\zeta
+\bar\psi_{\bar z}b \bar g_{\bar\zeta}+c.
\end{equation}
We solve (\ref{6}) together with its conjugate
as a system of two equations with the two unknowns
$g_{\bar\zeta}$ and $\bar g_\zeta$.
By ellipticity it has a unique solution.
In particular,
\begin{equation}
\label{27}
\tilde a=\frac{\psi_z}{\bar\psi_{\bar z}}\,
\frac{(a-\mu)(1-\bar a\mu)+|b|^2\mu}
{|1-\bar a\mu|^2-|b\mu|^2}.
\end{equation}
The equation $\tilde a=0$ turns into a quadratic equation
on $\mu$ of the form
\begin{equation}
\label{26}
\bar a\mu^2-(1+|a|^2-|b|^2)\mu+a=0.
\end{equation}
Due to $|a|+|b|<1$, the equation has two distinct solutions
$\mu_1, \mu_2$, $|\mu_1 \mu_2|=1$. We chose $\mu=\mu_1$,
the one with smaller modulus. (If $a(z)=0$ at some $z$,
then $\mu_1(z)=0$, $\mu_2(z)=\infty$.) It is easy to see
$\mu\in C^{k,\alpha}$ and $\|\mu\|_\infty<1$. Hence,
the homeomorphism $\psi\in C^{k+1,\alpha}$, the new
coefficients in (\ref{7}) are in $C^{k,\alpha}$,
and $\tilde a=0$.
$\blacksquare$
\medskip

Slightly changing notation, we now consider the equation
\begin{equation}
\label{8}
f_{\bar z}=A(z)\bar f_{\bar z}+b(z).
\end{equation}
Here for the sake of generality, $f$ and $b$
are $m$-vectors and $A$ is a $m\times m$ matrix, $m\ge1$.
To complete the proof of Theorem \ref{Dirichlet-scalar}
we need the following result only in the case
$\Omega=\D$, $f_0=0$.

\begin{thm}
\label{Dirichlet-chertov}
Let $\Omega\subset\C$ be a simply connected domain
of class $C^{k+1,\alpha}$, $k\ge0$, $0<\alpha<1$.
Let $A, b\in C^{k,\alpha}(\Omega)$,
$f_0\in C^{k+1,\alpha}(b\Omega)$,
$\|A\|_\infty<1$.
Then the equation (\ref{8}) with boundary condition
$\Re f|_{b\Omega}=f_0$ has a unique solution
in the Sobolev class $W^{1,2}(\Omega)$.
This solution $f\in C^{k+1, \alpha}(\Omega)$, and for fixed
$A$ the map $(b,f_0)\mapsto f$ is a bounded operator
$C^{k,\alpha}(\Omega)\times C^{k+1,\alpha}(b\Omega)\to
C^{k+1,\alpha}(\Omega)$.
\end{thm}

\proof
By Lemma \ref{lemma1} it suffices to
prove the result for $\Omega=\D$ and $f_0=0$.
For $f\in W^{1,2}(\D)$, the equation (\ref{8}) with boundary
conditions $\Re f|_{b\Omega}=0$ is equivalent to
\begin{equation}
\label{13}
f=T_1(A\bar f_{\bar z}+b).
\end{equation}
If this equation has a solution in $W^{1,2}(\D)$,
then $u=f_z$ satisfies the equation
\[
u=S_1(A\bar u+b).
\]
The latter by Proposition \ref{IntEquation-S1} has a unique
solution $u\in L^2(\D)$. This solution is in $C^{k,\alpha}(\D)$.
Then $f:=T_1(A\bar u+b)\in C^{k+1,\alpha}(\D)$ satisfies
(\ref{13}) because $f_z=S_1(A\bar u+b)=u$.
$\blacksquare$

\section{Inverting $f\mapsto f-T(af_z+b\bar f_{\bar z})$
in Lipschitz spaces}

Let $\Omega\subset\C$ be a bounded domain. We consider the integral
equation
\[
f=T(af_z+b\bar f_{\bar z})+c,
\]
here $a,b$, and $c$ are given functions in $\Omega$,
$f$ is unknown, and $T=T_\Omega$. Solving this equation
may be regarded as inverting the operator
$f\mapsto f-T(af_z+b\bar f_{\bar z})$.
Note that the solution
satisfies the boundary condition $Kf=Kc$ because $KT=0$,
here $K=K_\Omega$ is the Cauchy type integral.
To make the problem look similar to the one in the previous
section, we again consider the equation
\begin{equation}
\label{14}
f_{\bar z}=a(z)f_z+b(z)\bar f_{\bar z}+c(z)
\end{equation}
with boundary condition $Kf=Kf_0$ for a given
function $f_0$ on $b\Omega$. Our main result in the scalar
case is the following.

\begin{thm}
\label{InvertingT-scalar}
Let $\Omega\subset\C$ be a domain
of class $C^{k+1,\alpha}$, $k\ge0$, $0<\alpha<1$.
Let $a, b, c\in C^{k,\alpha}(\Omega)$,
$f_0\in C^{k+1,\alpha}(b\Omega)$,
$|a|+|b|\le a_0<1$, for some constant $a_0$.
Then the scalar equation (\ref{14}) with boundary condition
$Kf=Kf_0$ has a unique solution
in the Sobolev class $W^{1,2}(\Omega)$.
This solution $f\in C^{k+1,\alpha}(\Omega)$, and for fixed
$a$ and $b$, the map
$(c,f_0)\mapsto f$ is a bounded operator
$C^{k,\alpha}(\Omega)\times C^{k+1,\alpha}(b\Omega)\to
C^{k+1,\alpha}(\Omega)$.
\end{thm}

We again begin with reductions.
Note that $\Omega$ need not be simply connected, so instead
of the unit disc, we reduce to a domain of class $C^\infty$.
\begin{lemma}
\label{lemma5}
It suffices to prove Theorem \ref{InvertingT-scalar}
for $f_0=0$ and $\Omega$ of class $C^\infty$.
\end{lemma}
\proof
To deduce to $Kf=0$, we replace $f$ by $f-Kf_0$.
To reduce to a $C^\infty$-smooth domain $\Omega_0$
we again introduce a new independent
variable $\zeta=\psi(z)$ by a $C^{k+1,\alpha}$ diffeomorphism
$\psi:\bar\Omega\to\bar\Omega_0$.
To preserve the boundary condition $Kf=0$, we first
choose a conformal map
$\psi:\C\setminus\bar\Omega\to\C\setminus\bar\Omega_0$, and then
extend it $C^{k+1,\alpha}$-smoothly to $\Omega$.
(This procedure, however, will not preserve the conditions
$a=0$ or $b=0$ if they take place.)
$\blacksquare$
\medskip

In contrast to the proof of Theorem \ref{Dirichlet-scalar},
the reduction to $a=0$ is not straightforward because
a Beltrami homeomorphism \eqref{5} does not preserve the
boundary condition. Furthermore, the derivative of $\psi$
enters the boundary condition resulting in a loss of one derivative. We reduce to the case, in which $a$ is small.
\begin{lemma}
\label{lemma7} Under assumptions of Theorem \ref{InvertingT-scalar},
let $f_0=0$ and let $\Omega$ be $C^\infty$-smooth.
Let $\epsilon>0$.
There exists a $C^\infty$ diffeomorphism
$\psi:\C\to\C$ that transforms the equation \eqref{14}
in $\Omega$ into the equation \eqref{7} in $\Omega_0=\psi(\Omega)$,
in which $\|\tilde a\|_{C^{k,\alpha}}<\epsilon$.
The boundary condition $K_\Omega f=0$ transforms into
$K_{\Omega_0}g=Jg$, here $g=f\circ \psi^{-1}$ and
$J:L^2(b\Omega_0)\to C^\infty(b\Omega_0)$ is a smoothing
operator. Moreover, $\|\psi\|_{C^{k+1,\alpha}}$ is
bounded by a constant independent of $\epsilon$.
\end{lemma}
\proof
Let $C^\infty$-smooth functions $a_0$ and $b_0$ be
close to $a_0$ and $b_0$ in $C^{k,\alpha}(\Omega)$.
We find $\mu_0$ by solving \eqref{26} using $a_0$
and $b_0$ instead of $a$ and $b$.
We assume that $\mu_0$ is extended to the whole plane.
Following the proof of Lemma \ref{lemma3}
we make a substitution by a global Beltrami homeomorphism
of the equation \eqref{5} with $\mu_0$ instead of $\mu$.
If $a_0$ and $b_0$ are sufficiently close to
$a$ and $b$, then by \eqref{27}
the coefficient $\tilde a$ in \eqref{7},
satisfies  $\|\tilde a\|_{C^{k,\alpha}}<\epsilon$.

We now find out how the substitution affects the boundary
condition $K_\Omega f=0$.
Let $\rho$ be a defining function of $b\Omega_0$ with
$d\rho\neq 0$ in a neighborhood of $b\Omega_0$. Then
on $b\Omega_0$ we have
$\rho_\zeta d\zeta+\rho_{\bar\zeta}d\bar\zeta=0$.
Let $\zeta_0\in\Omega_0$ be sufficiently close to $b\Omega_0$.
With some abuse of notation we write
$z(\zeta)=\psi^{-1}(\zeta)$, $z=z(\zeta)$, $z_0=z(\zeta_0)$,
etc.
We have
\[
0=K_\Omega f(z_0)
=\frac{1}{2\pi i}\int_{b\Omega}\frac{f(z)\,dz}{z-z_0}
=\frac{1}{2\pi i}\int_{b\Omega_0}
\frac{g(\zeta)(z_\zeta d\zeta+z_{\bar\zeta}d\bar\zeta)}
{z-z_0}
=K_{\Omega_0}g(\zeta_0)-Jg(\zeta_0),
\]
\[
Jg(\zeta_0)=\frac{1}{2\pi i}\int_{b\Omega_0}
\(
\frac{1}{\zeta-\zeta_0}
-\frac{z_\zeta-\rho_{\bar\zeta}^{-1}\rho_\zeta z_{\bar\zeta}}
{z-z_0}
\)
g(\zeta)\,d\zeta.
\]
To understand the last integral we introduce
\begin{align*}
&\Phi(\zeta,\zeta_0)=z(\zeta)-z(\zeta_0)
-z_\zeta(\zeta)(\zeta-\zeta_0)
-z_{\bar\zeta}(\zeta)\bar{(\zeta-\zeta_0)},\\
&\phi(\zeta,\zeta_0)=\rho(\zeta)-\rho(\zeta_0)
-\rho_\zeta(\zeta)(\zeta-\zeta_0)
-\rho_{\bar\zeta}(\zeta)\bar{(\zeta-\zeta_0)}.
\end{align*}
We will use the last formula for $\zeta\in b\Omega_0$,
so we will have $\rho(\zeta)=0$. Then
$J=J_1+J_2$, here
\begin{align*}
&J_1g(\zeta_0)=\frac{1}{2\pi i}\int_{b\Omega_0}
\frac{\Phi(\zeta,\zeta_0)
-\rho_{\bar\zeta}^{-1}z_{\bar\zeta}\phi(\zeta,\zeta_0)}
{(\zeta-\zeta_0)(z-z_0)} g(\zeta)\,d\zeta\\
&J_2g(\zeta_0)=\frac{-\rho(\zeta_0)}{2\pi i}\int_{b\Omega_0}
\frac{
\rho_{\bar\zeta}^{-1}z_{\bar\zeta}g(\zeta)\,d\zeta}
{(\zeta-\zeta_0)(z-z_0)}.
\end{align*}
Since $\Phi(\zeta,\zeta_0)=O(|\zeta-\zeta_0|^2)$ and
$\phi(\zeta,\zeta_0)=O(|\zeta-\zeta_0|^2)$,
the kernel of the integral $J_1$ is $C^\infty$-smooth.
Hence $J_1$ is a smoothing operator.
The integral $J_2$ reduces to the integral $K^0_2$
introduced in the proof of Lemma \ref{lemma6} by \eqref{31}.
Since $g\in C^{k,\alpha}(\Omega_0)\subset C^\alpha(\Omega_0)$,
the argument in the proof of that lemma implies
$J_2g(\zeta_0)=O(|\rho(\zeta_0)|^\alpha)$.
This estimate means that $J_2g$ has zero boundary
values on $b\Omega_0$, and $J$ reduces to the smoothing
operator $J_1$.

Finally,
the $C^{k+1,\alpha}$ norm of $\psi$ depends only
on $C^{k,\alpha}$ norm of $\mu_0$ (see \cite{AIM,V}), which in turn depends
only on $C^{k,\alpha}$ norms of $a$ and $b$, hence
the last assertion in the lemma will hold automatically.
The lemma is proved.
$\blacksquare$
\medskip

We again state a special case of Theorem \ref{InvertingT-scalar}
in a vector from.
Slightly changing notation, we now consider the equation
\begin{equation}
\label{28}
f_{\bar z}=A_1(z)f_z+A_2(z)\bar f_{\bar z}+b(z).
\end{equation}
Here $f$ and $b$ are $m$-vectors and $A_1$ and $A_2$ are
$m\times m$ matrix, $m\ge1$.
To complete the proof of Theorem \ref{InvertingT-scalar}
it suffices to prove the following result.

\begin{thm}
\label{InvertingT-chertov}
Let $\Omega\subset\C$ be a $C^\infty$-smooth domain.
Let $A_2, b\in C^{k,\alpha}(\Omega)$, $k\ge0$, $0<\alpha<1$.
Suppose $\|A_2\|_\infty<1$.
Then there exists $\epsilon>0$ such that if
$A_1\in{C^{k,\alpha}(\Omega)}$ and
$\|A_1\|_{C^{k,\alpha}(\Omega)}<\epsilon$, then
for every $f_0\in C^{k+1,\alpha}(b\Omega)$,
the equation (\ref{28}) with boundary condition
$K_\Omega f=K_\Omega f_0$ has a unique solution
in the Sobolev class $W^{1,2}(\Omega)$.
This solution $f\in C^{k+1, \alpha}(\Omega)$, and for fixed
$A_1$ and $A_2$, the map
$(b,f_0)\mapsto f$ is a bounded operator
$C^{k,\alpha}(\Omega)\times C^{k+1,\alpha}(b\Omega)\to
C^{k+1,\alpha}(\Omega)$.
\end{thm}
\proof
By Lemma \ref{lemma5} we assume $f_0=0$.
Let $T=T_\Omega$ and $S=S_\Omega$.
For $f\in W^{1,2}(\Omega)$, the equation (\ref{28})
with boundary conditions $K_\Omega f=0$ is equivalent to
\begin{equation}
\label{29}
f=T(A_1(z)f_z+A_2(z)\bar f_{\bar z}+b).
\end{equation}
If this equation has a solution in $W^{1,2}(\Omega)$,
then $u=f_z$ satisfies the equation
\begin{equation}
\label{30}
u=S(A_1u+A_2\bar u+b).
\end{equation}
If $A_1$ is small, then of course \eqref{30}
has a unique solution $u\in L^2(\Omega)$.
By the proof of Proposition \ref{IntEquation-SD}
the equation \eqref{30} implies
\[
u=Mu+NS(A_1u+b).
\]
Since $NS$ is bounded in $C^{k,\alpha}(\Omega)$,
we choose $0<\epsilon< \|NS\|_{C^{k,\alpha}(\Omega)}^{-1}$.
Then if
$\|A_1\|_{C^{k,\alpha}(\Omega)}<\epsilon$,
then $(I-NSA_1)^{-1}$ is bounded in $C^{k,\alpha}(\Omega)$
and $L^2(\Omega)$, and
\[
u=M_1u+N_1b,\qquad
M_1=(I-NSA_1)^{-1}M,\qquad
N_1=(I-NSA_1)^{-1}NS.
\]
By bootstrapping $u\in C^{k,\alpha}(\Omega)$.
Then $f:=T(A_1u+A_2\bar u+b)\in C^{k+1,\alpha}(\Omega)$ satisfies
(\ref{28}) because $f_z=S(A_1u+A_2\bar u+b)=u$.
$\blacksquare$
\medskip

\noindent
{\bf Proof of Theorem \ref{InvertingT-scalar}.}
There is a unique solution $f\in W^{1,2}(\Omega)$
of \eqref{14} with $K_\Omega f=0$.
Indeed, as we argued before, since $|a|+|b|\le a_0<1$,
there is a unique
$u\in L^2(\Omega)$ satisfying $u=S(au+b\bar u+c)$.
Then $f=T(au+b\bar u+c)$. By Lemma \ref{lemma7},
after the substitution $\zeta=\psi(z)$, the function
$g=f\circ\psi^{-1}$ satisfies \eqref{7} with small
$\|\tilde a\|_{C^{k+1,\alpha}}$ and $K_{\Omega_0}g=Jg$.
Since $K_\Omega f=0$, of course $Jg$ is holomorphic.
Since $Jg$ is $C^\infty$, by Theorem \ref{InvertingT-chertov}
we have $g\in C^{k+1,\alpha}(\Omega_0)$.
Hence $f\in C^{k+1,\alpha}(\Omega)$.
$\blacksquare$
\medskip

In conclusion we point out that Theorem \ref{InvertingT-chertov}
answers a question raised in \cite{BGR}.
Let $A$ be a $m\times m$ matrix function of class
$C^{k,\alpha}(\Omega)$, $k\ge0$, $0<\alpha<1$,
in a $C^\infty$-smooth bounded
domain $\Omega\subset\C$, $\|A\|_\infty<1$, $m\ge2$.
The question from \cite{BGR} (Problem B)
reduces to asking whether the operator
$f\mapsto f-T_\Omega (A \bar{\partial f})$
has a bounded inverse in $C^{k+1,\alpha}(\Omega)$.
The affirmative answer
is given by Theorem \ref{InvertingT-chertov}
with $A_1=0$, $A_2=A$. The authors also raise
a similar question (Problem A) for the operator
$f\mapsto f-T_\Omega (A {\partial f})$.
However, in this paper we are able to treat this
question only in the scalar case
(Theorem \ref{InvertingT-scalar}).


\begin{thebibliography}{CIT}

\bibitem{ADN}
S. Agmon, A. Douglis and L. Nirenberg,
{\it Estimates near the boundary for solutions of elliptic partial
differential equations satisfying general boundary conditions. I},
Comm. Pure Appl. Math. 12 (1959), 623--727.

\bibitem{AIM}
K. Astala, T. Iwaniec, and G. Martin,
{\it Elliptic partial differential equations
and quasiconformal mappings in the plane},
Princeton Mathematical Series, 48, Princeton, 2009, xviii+677 pp.

\bibitem{Bell}
S. Bell,
{\it The Cauchy transform, potential theory, and conformal mapping},
CRC Press, Boca Raton, FL, 1992. x+149 pp.

\bibitem{BGR}
F. Bertrand, X. Gong, and J.-P. Rosay,
{\it Common boundary values of holomorphic
functions for two-sided complex structures},
arXiv 1008.1234.

\bibitem{C}
A.-P. Calder\'on,
{\it Commutators of singular integral operators},
Proc. Nat. Acad. Sci. U.S.A. {\bf 53} (1965), 1092--1099.

\bibitem{CRW}
R. R. Coifman, R. Rochberg, and G. Weiss,
{\it Factorization theorems for Hardy spaces
in several variables},
Ann. of Math. (2) {\bf 103} (1976), 611--635.

\bibitem{IvRo}
S. Ivashkovich and J.-P. Rosay,
{\it Schwarz-type lemmas for $\bar\partial$-inequalities and
complete hyperbolicity of almost complex manifolds},
Ann. Inst. Fourier {\bf 54} (2004), 2387--2435.

\bibitem{KS}
N. Kerzman and E. M. Stein,
{\it The Cauchy kernel, the Szeg\"o kernel,
and the Riemann mapping function},
Math. Ann. 236 (1978), 85--93.

\bibitem{NiWo}
A. Nijenhuis and W. Woolf,
{\it Some integration problems in almost-complex
and complex manifolds},
Ann. of Math. (2) {\bf 77} (1963), 429--484.

\bibitem{P}
S. Pr\"ossdorf,
{\it Some classes of singular equations},
North-Holland Publishing Co.,
Amsterdam-New York, 1978. xiv+417 pp.

\bibitem{S}
J. Schauder,
{\it \"Uber lineare elliptische Differentialgleichungen
zweiter Ordnung},
Math. Z. 38 (1934), 257--282.

\bibitem{ST}
A. Sukhov and A. Tumanov,
{\it Boundary value problems and J-complex curves},
Complex Var. Elliptic Equ. 58 (2013), 1549--1557.

\bibitem{V}
I. N. Vekua, {\it Generalized analytic functions},
Moscow, 1959; English translation: Pergamon Press, London, 1962.


\end{thebibliography}
\end{document}